\documentclass[10pt,a4paper]{amsart}
\usepackage{amssymb}%,showkeys
\newtheorem{theorem}{Theorem}[section]
\newtheorem{proposition}[theorem]{Proposition}
\newtheorem{lemma}[theorem]{Lemma}

\theoremstyle{remark}
\newtheorem*{remark}{Remark}
\newtheorem*{remarks}{Remarks}

\newtheorem*{definition}{Definition}

\newcommand{\remove}[1]{ }
\newcommand{\set}[1]{\left\{#1\right\}}
\def\bfu{\mathbf U}
\def\bfuu{\overline{\mathbf U}}
\def\bfv{\mathbf V}
\def\bfj{\mathbf J}
\def\bfw{\mathbf W}
\def\ff{\mathcal F}
\def\hh{\mathcal H}
\def\uu{\mathcal U}
\def\vv{\mathcal V}
\def\uuu{\overline{\mathcal U}}
\def\uuuq{\overline{\mathcal{U}_q}}
\def\NN{\mathbb N}

\def\RR{\mathbb R}

\numberwithin{equation}{section}

\begin{document}

\title{A two-dimensional univoque set}
\author{Martijn de Vries}
\address{Delft University of Technology, Mekelweg 4, 2628 CD Delft, the Netherlands}
\email{w.m.devries@tudelft.nl}
\author{Vilmos Komornik}
\address{D\'epartement de Math\'ematique,
         Universit\'e de Strasbourg,
         7 rue Ren\'e Descartes, 67084 Strasbourg Cedex, France}
\email{komornik@math.u-strasbg.fr}
\subjclass[2000]{Primary:11A63, Secondary:11B83}
\keywords{Greedy expansion, beta-expansion,
univoque sequence, univoque set, Cantor set, Hausdorff dimension.}
\date{\today}

\begin{abstract}
Let $\bfj \subset \RR^2$ be the set of couples $(x,q)$ with $q>1$ such that $x$ has at least one representation of the form $x=\sum_{i=1}^{\infty} c_i q^{-i}$ with integer coefficients $c_i$ satisfying $0 \le c_i < q$, $i \ge 1$. In this case we say that $(c_i)=c_1c_2\ldots$ is an expansion of $x$ in base $q$. Let $\bfu$ be the set of couples $(x,q) \in \bf J$ such that $x$ has exactly one expansion in base $q$.

In this paper we deduce some topological and combinatorial properties of the set $\bfu$.
We characterize the closure of $\bfu$, and we determine its Hausdorff dimension. For $(x,q) \in \bfj$, we also prove new properties of the lexicographically largest expansion of $x$ in base $q$.
\end{abstract}

\maketitle

\section{Introduction}\label{s1}

Let $\bfj$ be the set consisting of all elements $(x,q)\in \mathbb{R} \times (1, \infty)$ such that there exists at least one sequence $(c_i)=c_1c_2\ldots$ of integers satisfying $0\le c_i<q$ for all $i$, and
\begin{equation}\label{11}
x=\frac{c_1}{q} + \frac{c_2}{q^2}+\cdots.
\end{equation}
If \eqref{11} holds, we say that $(c_i)$ is an \emph{expansion of $x$ in base $q$}, and if the base $q$ is understood from the context, we sometimes simply say that $(c_i)$ is an expansion of $x$. The numbers $c_i$ of an expansion $(c_i)$ are usually referred to as {\it digits}. We denote by $\lceil q \rceil$ the smallest integer larger than or equal to $q$. The {\it alphabet} $A_q$ is the set of ``admissible'' digits in base $q$, i.e.,
$A_q= \set{0, \ldots, \lceil q \rceil - 1}$.

If $q>1$ and $0 \le x \le (\lceil q \rceil - 1)/(q-1)$, then a particular expansion of $x$ in base $q$, the so-called \emph{quasi-greedy expansion} $(a_i(x,q))$, may be defined recursively as follows.
For $x=0$ we set $(a_i(x,q)):=0^{\infty}$. If $x>0$ and $a_i(x,q)$ has already been defined for
$1 \le i<n$ (no condition if $n=1$), then $a_n(x,q)$ is the largest element of $A_q$ satisfying
\begin{equation*}
\frac{a_1(x,q)}{q} + \cdots + \frac{a_n(x,q)}{q^n} <x.
\end{equation*}
One easily verifies that $(a_i(x,q))$ is indeed an expansion of $x$ in base $q$. Therefore
\begin{equation*}
(x,q)\in\bfj \Longleftrightarrow q>1\quad\text{and}\quad  x \in J_q:= \left[0, \frac{\lceil q\rceil -1}{q-1}\right].
\end{equation*}

Let us denote by $\bfu $ the set of couples $(x,q) \in \bfj$ such that $x$ has \emph{exactly one} expansion in base $q$. For
example, $(0,q)\in\bfu $ for every $q> 1$, but $\bfu $ has much more elements. The main purpose of this paper is to
describe the topological and combinatorial nature of $\bfu$. We will prove the following theorem:

\begin{theorem}\label{t11}\mbox{}

\begin{itemize}
\item[\rm (i)] The set $\bfu $ is not closed. Its closure $\bfuu $ is a Cantor set \footnote{We recall that a Cantor set is a nonempty closed set having neither interior nor isolated points.}.
\item[\rm (ii)] Both $\bfu $ and $\bfuu $ are two-dimensional Lebesgue null sets.
\item[\rm (iii)] Both $\bfu $ and $\bfuu $  have Hausdorff dimension two.
\end{itemize}
\end{theorem}

As far as we know this two-dimensional {\it univoque} set has not yet been investigated.
There exists, however, a number of papers devoted to the study of its one-dimensional sections
\begin{equation*}
\uu:=\set{q>1\ :\ (1,q)\in\bfu}
\end{equation*}
and
\begin{equation*}
\uu_q:=\set{x \in J_q :\ (x,q)\in\bfu },\quad q>1.
\end{equation*}

The study of $\uu$ started with the paper of Erd\H os, Horv\'ath and Jo\'o \cite{EHJ} and was studied subsequently in \cite{DK2}, \cite{DVK}, \cite{EJK1}, \cite{EJK2}, \cite{KK}, \cite{KL1}, \cite{KL2}. We recall
in particular that $\uu$ and its closure $\uuu$ have Lebesgue measure zero and Hausdorff dimension one.

The sets $\uu_q$
have been investigated in \cite{DK1}, \cite{DK2}, \cite{DVK}, \cite{GS}, \cite{K1}, \cite{K2}.
It is known that $\uu_q$ is closed if and only if $q$ does not belong to the null set $\uuu$, and
that its closure $\uuuq$ has Lebesgue measure zero for all non-integer bases $q>1$. Moreover, the set of numbers $x \in J_q$ having continuum many expansions in base $q$ has full Lebesgue measure for each non-integer $q>1$  (see \cite{DDV}, \cite{Si1}, \cite{Si2}).

The key to the proof of Theorem~\ref{t11} is an algebraic characterization of $\bfuu$ by using the quasi-greedy expansions
$(a_i(x,q))$. We write for brevity $\alpha_i(q):=a_i(1,q)$, $i \in \NN:=\set{1,2,\ldots}, q>1$. Note that $\alpha_1(q)=\lceil q\rceil -1$, the largest admissible digit in base $q$. In the statement of the following theorem we use the lexicographic
order between sequences and we define the {\it conjugate} of the number $a_i(x,q)$ by  $\overline{a_i(x,q)}:=\alpha_1(q) - a_i(x,q)$. If $q>1$ and $c_i \in A_q$, $ i \ge 1$, we shall also write $\overline{c_1 \ldots c_n}$ instead of $\overline{c_1} \ldots \overline{c_n}$ and $\overline{c_1c_2 \ldots}$ instead of $\overline{c_1} \, \overline{c_2} \ldots$.

\begin{theorem}\label{t12}
A point $(x,q)\in\bfj $ belongs to $\bfuu $ if and only if
\begin{equation*}
\overline{a_{n+1}(x,q)a_{n+2}(x,q)\ldots }\le \alpha_1(q)\alpha_2(q)\ldots \quad\text{whenever }a_n(x,q)>0.
\end{equation*}
\end{theorem}

Along with the quasi-greedy expansion, we also need the notion of the \emph{greedy expansion} $(b_i(x,q))$ for $x \in J_q$, introduced by R\'enyi \cite{R}. It can be defined by a slight modification of the  above recursion: if $b_i(x,q)$ has already been defined for all $1 \le i<n$ (no condition if $n=1$), then $b_n(x,q)$ is the largest element of $A_q$  satisfying
\begin{equation*}
\frac{b_1(x,q)}{q}+ \cdots +  \frac{b_n(x,q)}{q^n} \le x.
\end{equation*}
Note that the  greedy expansion  $(b_i(x,q))$ of a number $x \in J_q$ is the lexicographically largest expansion of $x$ in base $q$. We denote the greedy expansion of $1$ in base $q$ by $(\beta_i(q)):= (b_i(1,q))$.

The rest of this paper is organized as follows. In the next section we give a short overview of some basic results on greedy and quasi-greedy expansions, and we prove some new results concerning the coordinate-wise convergence of sequences of these expansions. We shall prove (see Theorem~\ref{t27}) that the set of numbers $x \in J_q$ for which the greedy expansion of $x$ in base $q$ is not the greedy expansion of a number belonging to $J_p$ in any smaller  base $p \in (1,q)$ is of full Lebesgue measure and its complement in $J_q$ is a set of first category and  Hausdorff dimension one. We shall also prove (see Theorem~\ref{t28}) that for each word $v:=b_{\ell+1}(x,q) \ldots b_{\ell+m}(x,q)$ $(\ell \ge 0, m \ge 1, x \in [0,1))$ there exists a set $Y_v \subset J_q$ of first category and Hausdorff dimension less than one, such that the word $v$ occurs in the greedy expansion in base $q$ of every number belonging to $J_q \setminus Y_v$. Using (some of) the results of Section~\ref{s2} we prove Theorem~\ref{t12} in Section~\ref{s3} and Theorem~\ref{t11} in Section~\ref{s4}.

\section{Greedy and quasi-greedy expansions}\label{s2}

In this paper we consider only one-sided sequences of nonnegative integers. We equip this set of sequences $\set{0,1, \ldots}^{\NN}$ with the topology of coordinate-wise convergence. We say that an expansion is \emph{infinite} if it has infinitely many nonzero elements; otherwise it is called {\it finite}. Using this terminology, the quasi-greedy expansion $(a_i(x,q))$ of a number $x \in J_q \setminus \set{0}$
is the lexicographically largest \emph{infinite} expansion of $x$ in base $q$.

The family of all quasi-greedy expansions is characterized by the following propositions (see \cite{BK} or \cite{DVK} for a proof):

\begin{proposition}\label{p21}
The map $q \mapsto (\alpha_i(q))$ is a strictly increasing bijection from
the open interval $(1, \infty)$ onto the set of all {\rm infinite} sequences $(\alpha_i)$
satisfying
\begin{equation*}
\alpha_{k+1} \alpha_{k+2} \ldots \leq \alpha_1 \alpha_2 \ldots \quad
\mbox{for all} \quad k\geq 1.
\end{equation*}
\end{proposition}

\begin{proposition}\label{p22}
For each $q>1$, the map $x \mapsto (a_i(x,q))$ is a strictly increasing bijection from $J_q \setminus \set{0}$ onto the set of all {\rm infinite} sequences $(a_i)$ satisfying
\begin{equation*}
a_n \in A_q \quad \text{for all} \quad n \ge 1
\end{equation*}
and
\begin{equation*}
a_{n+1}a_{n+2}\ldots \le \alpha_1(q)\alpha_2(q)\ldots \quad\text{whenever }a_n <\alpha_1(q).
\end{equation*}
\end{proposition}

The quasi-greedy expansions have a lower semicontinuity property for the order topology induced by the lexicographic
order. More precisely, we have the following result.

\begin{lemma}\label{l23}
Let $(x,q)\in \bfj $. Then

\begin{itemize}
\item[\rm (i)] for each positive integer $m$
there exists a neighborhood $\bfw \subset \RR^2$ of $(x,q)$ such that
\begin{equation}\label{21}
a_1(y,r)\ldots a_m(y,r)\ge a_1(x,q)\ldots a_m(x,q) \quad \text{for all} \quad  (y,r)\in\bfw \cap \bfj;
\end{equation}
\item[\rm (ii)] if  $(y_n,r_n)$ converges to $(x,q)$ in $\bfj $ from below, then
$(a_i(y_n,r_n))$ converges to $(a_i(x,q))$.
\end{itemize}
\end{lemma}

\begin{proof}
(i) We may assume that $x \not=0$. By definition of the quasi-greedy expansion we have
\begin{equation*}
\sum_{i=1}^n\frac{a_i(x,q)}{q^i}<x \text{ for all } n=1,2,\ldots.
\end{equation*}
For any fixed positive integer $m$, if $(y,r)\in\bfj $ is sufficiently close to $(x,q)$, then
$r > \lceil q \rceil -1$, i.e., $A_q \subset A_r$, and
\begin{equation*}
\sum_{i=1}^n\frac{a_i(x,q)}{r^i}<y,\quad n=1,\ldots, m.
\end{equation*}
These inequalities imply \eqref{21}.

(ii) If $y_n\le x$ and $r_n\le q$, we deduce from the definition of the quasi-greedy expansion that
\begin{equation*}
(a_i(x,q))\ge (a_i(y_n,r_n))
\end{equation*}
for every $n$. Equivalently, we have
\begin{equation*}
a_1(x,q)\ldots a_m(x,q)\ge a_1(y_n,r_n)\ldots a_m(y_n,r_n)
\end{equation*}
for all positive integers $m$ and $n$. It remains to notice that by the previous part the converse inequality also
holds for each fixed $m$ if $n$ is large enough.
\end{proof}

\remove{
\begin{remarks}
We complete the preceding remarks. \mbox{}

\begin{itemize}
\item If  $(y_n,r_n)$ converges to $(x,q)$ in $\bfj $ from below and if $(y_n,r_n)\ne (x,q)$ for every $n$, then
$(b_i(y_n,r_n))$ also converges to $(a_i(x,q))$. This follows from the first remark above by observing that
for $(y_n,r_n)\ne (x,q)$ the inequalities
\begin{equation*}
(a_i(x,q))\ge (b_i(y_n,r_n))\ge (a_i(y_n,r_n))
\end{equation*}
hold by the definition of quasi-greedy and greedy expansions.

\item Applying the preceding remark with $x=y_n=1$ we obtain that if $r_n\uparrow q$ in $(1,\infty)$, and $r_n\ne q$ for every
$n$, then $(\beta_i(r_n))$ also converges to $(\alpha_i(q))$.
\end{itemize}
\end{remarks}
}

The family of greedy expansions has already been characterized by Parry \cite{P}:

\begin{proposition}\label{p24}
For a given base $q>1$, the map $x \mapsto (b_i(x,q))$ is a strictly increasing bijection from $J_q$ onto the set of all sequences $(b_i)$ satisfying
\begin{equation*}
b_n \in A_q \quad \text{for all} \quad n \ge 1
\end{equation*}
and
\begin{equation*}
b_{n+1}b_{n+2}\ldots <\alpha_1(q)\alpha_2(q)\ldots \quad \text{whenever} \quad b_n< \alpha_1(q).
\end{equation*}
\end{proposition}

The greedy expansions have the following upper semicontinuity property:

\begin{lemma}\label{l25}
Let $(x,q)\in \bfj $ and suppose $q$ is a non-integer. Then
\begin{itemize}
\item[\rm (i)] for each positive integer $m$
there exists a neighborhood $\bfw \subset \RR^2$ of $(x,q)$ such that
\begin{equation}\label{22}
b_1(y,r)\ldots b_m(y,r)\le b_1(x,q)\ldots b_m(x,q) \quad \text{for all} \quad (y,r)\in\bfw \cap \bfj;
\end{equation}
\item[\rm (ii)] if  $(y_n,r_n)$ converges to $(x,q)$ in $\bfj $ from above, then
$(b_i(y_n,r_n))$ converges to $(b_i(x,q))$.
\end{itemize}
\end{lemma}

\begin{proof}
(i) By the definition of greedy expansions we have
\begin{equation*}
\sum_{i=1}^n\frac{b_i(x,q)}{q^i}>x-\frac{1}{q^n}\quad\text{whenever}\quad b_n(x,q)< \alpha_1(q).
\end{equation*}
If $(y,r)\in\bfj $ is sufficiently close to $(x,q)$, then $A_q = A_r$, $\alpha_1(r)=\alpha_1(q)$, and
\begin{equation*}
\sum_{i=1}^n\frac{b_i(x,q)}{r^i}>y-\frac{1}{r^n}\quad\text{whenever $n\le m$ and}\quad b_n(x,q)< \alpha_1(r).
\end{equation*}
These inequalities imply \eqref{22}.

(ii) If $y_n\ge x$ and $r_n\ge q$, we deduce from the definition of the greedy expansion that
\begin{equation*}
(b_i(x,q))\le (b_i(y_n,r_n))
\end{equation*}
for every $n$. Equivalently, we have
\begin{equation*}
b_1(x,q)\ldots b_m(x,q)\le b_1(y_n,r_n)\ldots b_m(y_n,r_n)
\end{equation*}
for all positive integers $m$ and $n$. It remains to notice that by the previous part the converse inequality also
holds for each fixed $m$ if $n$ is large enough.
\end{proof}

\remove{
\begin{remarks}\mbox{}

If, moreover, $(y_n,r_n)\ne (x,q)$ for every $n$, then $(a_i(y_n,r_n))$ also converges to
$(b_i(x,q))$. This follows from the preceding remark by observing that for $(y_n,r_n)\ne (x,q)$ the inequalities
\begin{equation*}
(b_i(x,q))\le (a_i(y_n,r_n))\le (b_i(y_n,r_n))
\end{equation*}
also hold by the definition of quasi-greedy expansions.

Applying the first remark with $x=y_n=1$ we get that if  $r_n\downarrow q$ in $(1,\infty)$, then $(\beta_i(r_n))$
converges to $(\beta_i(q))$.

Applying the second remark with $x=y_n=1$ we get that if, moreover, $r_n\ne q$ for every $n$, then $(\alpha_i(r_n))$
also converges to $(\beta_i(q))$.
\end{remarks}
}
From Lemmas \ref{l23} and \ref{l25} we deduce the following result:

\begin{proposition}\label{p26}
Consider $(x,q)\in \bfj $ with a non-integer base $q$ and assume that the greedy expansion $(b_i(x,q))$ is infinite.
If $(y_n,r_n)$ converges to $(x,q)$ in $\bfj$, then both
$(a_i(y_n,r_n))$ and $(b_i(y_n,r_n))$ converge to $(b_i(x,q))=(a_i(x,q))$.
\end{proposition}
\begin{proof}
For each positive integer $m$ there exists a neighborhood $\bfw \subset \RR^2$ of $(x,q)$ such that for all $(y,r)\in\bfw \cap \bfj$,
\begin{align*}
a_1(x,q)\ldots a_m(x,q) &\le a_1(y,r)\ldots a_m(y,r) \\
& \le b_1(y,r)\ldots b_m(y,r) \\
& \le b_1(x,q)\ldots b_m(x,q).
\end{align*}
The result follows from our assumption that $(a_i(x,q))=(b_i(x,q))$.
\end{proof}

\begin{theorem}\label{t27}
Let $q>1$ be a real number. Then
\begin{itemize}
\item[\rm (i)] for each $r \in (1,q)$, the Hausdorff dimension of the set
\begin{equation*}
G_{r,q}:=\set{\sum_{i=1}^{\infty} \frac{b_i(x,r)}{q^i} : x \in J_r}
\end{equation*}
equals $\log r / \log q$;
\item[\rm (ii)] the set
\begin{equation*}
G_q:=\bigcup \set{G_{r,q}: r \in (1,q)}
\end{equation*}
is of first category, has Lebesgue measure zero and Hausdorff dimension one.
\end{itemize}
\end{theorem}

\begin{proof} (i)
It is well known (\cite{KL2}, \cite{P}) and easy to prove that the set of numbers $r>1$ for which $(\beta_i(r))$ is finite is dense in $[1, \infty)$. Moreover, if $(\beta_i(r))$ is finite and $\beta_n(r)$ is its last nonzero element, then $(\alpha_i(r))= (\beta_1(r) \ldots \beta_{n-1}(r) \beta_n^-(r))^{\infty}$ $(\beta_n^-(r):=\beta_n(r)-1)$. By virtue of Propositions~\ref{p21} and \ref{p24} we have $G_{s,q} \subset G_{t,q}$ whenever $1 < s < t < q$. Hence it is enough to prove that ${\rm dim}_{H} G_{r,q} = \log r / \log q$ for those values $r \in (1,q)$ for which $(\alpha_i(r))$ is periodic.

Fix $r \in (1,q)$ such that $(\alpha_i):=(\alpha_i(r))$ is periodic and let $n \in \NN$ be such that $(\alpha_i)= (\alpha_1 \ldots \alpha_n)^{\infty}$. Let us denote by $W_r$ the set consisting of the finite words
\begin{equation*}
w_{ij}:=\alpha_1 \ldots \alpha_{j-1}i, \quad 0 \le i < \alpha_j, \quad 1 \le j \le n
\end{equation*}
and
\begin{equation*}
w_{\alpha_n n}:=\alpha_1 \ldots \alpha_{n-1}\alpha_n.
\end{equation*}
Let $\ff_r'$ be the set of sequences $(c_i)=c_1 c_2$ such that for each $k \ge 0$ the inequality $c_{k+1} \ldots c_{k+n} \le \alpha_1 \ldots \alpha_n$ holds. Note that the set $\ff_r'$ consists of those sequences $(c_i)$ such that each tail of $(c_i)$ (including $(c_i)$ itself) starts with a word belonging to $W_r$. It follows from Propositions~\ref{p21} and \ref{p24} that a sequence $(b_i)$ is greedy in base $r$ if and only if $b_m \in A_r$ for all $m \ge 1$ and
\begin{equation*}
b_{m+k+1}b_{m+k+2} \ldots < \alpha_1 \alpha_2 \ldots
\quad \text{for all $k \ge0$, whenever } b_m < \alpha_1.
\end{equation*}
Therefore, any greedy expansion $(b_i) \not= \alpha_1^{\infty}$ in base $r$ can be written as $\alpha_1^{\ell} c_1c_2 \ldots$ for some $\ell \ge 0$ ($\alpha_1^0$ denotes the empty word) and some sequence $(c_i)$ belonging to $\ff_r'$. Conversely, if no tail of a sequence belonging to $\ff_r'$  equals $(\alpha_i)$, then it is the greedy expansion in base $r$ of some $x \in J_r$. Hence if we set
\begin{equation*}
\ff_{r,q}:= \set{\sum_{i=1}^{\infty} \frac{c_i}{q^i} : (c_i) \in \ff_r'},
\end{equation*}
then $\ff_{r,q} \setminus G_{r,q}$ is countable and $G_{r,q}$ can be covered by countably many sets similar to $\ff_{r,q}$.
Since  the union of countably many sets of Hausdorff dimension $s$ is still of Hausdorff dimension $s$, we have ${\rm dim}_{H}G_{r,q}={\rm dim}_{H} \ff_{r,q}$.

We associate with each word $w_{ij}\in W_r$ a similarity $S_{ij}: J_q \to J_q$ defined by the formula
\begin{equation*}
S_{ij}(x):=\frac{\alpha_1}{q}+\cdots+\frac{\alpha_{j-1}}{q^{j-1}}+\frac{i}{q^j}+\frac{x}{q^j}, \quad x \in J_q.
\end{equation*}
It follows from Proposition~\ref{p21} and the definition of $\ff_{r,q}$ that
\begin{equation}\label{23}
\ff_{r,q}=\bigcup S_{ij}(\ff_{r,q})
\end{equation}
where the union runs over all $i$ and $j$ for which $w_{ij} \in W_r$. Applying Proposition~\ref{p21} again, it follows that $r$ is the largest element of the set of numbers $t>1$ for which $\alpha_i(t)=\alpha_i$,  $1 \le i \le n$. Hence $\alpha_1 \ldots \alpha_n < \alpha_1(q) \ldots \alpha_n(q)$ and therefore each sequence in $\ff_r'$ is the greedy expansion in base $q$ of some $x \in \ff_{r,q}$. It follows that the sets $S_{ij}(\ff_{r,q})$ on the right side of \eqref{23} are disjoint. Moreover, the function $x \mapsto (b_i(x,q))$ that maps $\ff_{r,q}$ onto $\ff_r'$ is increasing. Using the definition of $\ff_r'$ it is easily seen that the limit of each monotonic sequence of elements in $\ff_{r,q}$ belongs to $\ff_{r,q}$. We conclude that the closed set $\ff_{r,q}$ is the (nonempty compact) invariant set of this system of similarities.
An application of Propositions 9.6 and 9.7 in \cite{F} yields that
\begin{equation*}
{\rm dim}_{H} \ff_{r,q} = {\rm dim}_{H} G_{r,q} = s
\end{equation*}
where $s$ is the real solution of the equation
\begin{equation*}
\frac{\alpha_1}{q^s}+\cdots+\frac{\alpha_{n-1}}{q^{(n-1)s}}+ \frac{\alpha_n+1}{q^{ns}}=1.
\end{equation*}
Since
\begin{equation*}
\frac{\alpha_1}{r} + \cdots + \frac{\alpha_{n-1}}{r^{n-1}} + \frac{\alpha_n + 1}{r^n}=1
\end{equation*}
 we have $s= \log r / \log q$.

(ii) It follows at once from Theorem~\ref{t27}(i) that ${\rm dim}_{H}G_q=1$. Let $r \in (1,q)$ be such that $(\alpha_i(r))$ is periodic. The proof of Theorem~\ref{t27}(i) shows that
\begin{equation*}
G_{r,q} \subset \bigcup_{n=1}^{\infty} \left(a_n + b_n \ff_{r,q}\right)
\end{equation*}
for some constants $a_n, b_n \in \RR$ $(n \in \NN)$. Since $\ff_{r,q}$ is a closed set of Hausdorff dimension less than one, it follows in particular that the sets $a_n + b_n \ff_{r,q}$ are nowhere dense null sets. Since $G_{s,q} \subset G_{t,q}$ whenever $1 < s < t < q$, the set $G_q$ is a null set of first category.
\end{proof}

\begin{theorem}\label{t28} Let $q>1$ be a real number.
\begin{itemize}
\item[\rm (i)] Let $v:=b_{\ell+1}(y,q) \ldots b_{\ell+m}(y,q)$ for some $y \in [0,1)$ and some integers $\ell \ge 0$ and $m \ge 1$. The set $Y_v$ of numbers $x \in J_q$ for which the word $v$ does not occur in the greedy expansion of $x$ in base $q$ has Hausdorff dimension less than one.

\item[\rm (ii)] The set $Y$ of numbers $x \in J_q$ for which {\rm at least one} word  of the form $b_{\ell+1}(y,q) \ldots b_{\ell+m}(y,q)$ $(\ell \ge 0$, $m \ge 1$, $y \in [0,1))$ does not occur in the greedy expansion of $x$ in base $q$ is of first category, has Lebesgue measure zero and Hausdorff dimension one.
\end{itemize}
\end{theorem}

\begin{proof}
(i) Using the inequality $(b_i(y,q)) < (\alpha_i(q))$, it follows from Proposition~\ref{p24} that for some $k \in \NN$, there exist positive integers $m_1, \ldots,m_k$ and nonnegative integers $\ell_1, \ldots, \ell_k$ satisfying $\alpha_{m_j}(q)>0$ and $\ell_j < \alpha_{m_j}(q)$ for each $1 \le j \le k$, such that $v$ is a subword of
\begin{equation*}
w:= \alpha_1(q) \ldots \alpha_{m_1-1}(q)\ell_1 \ldots \alpha_1(q) \ldots \alpha_{m_k-1}(q)\ell_k.
\end{equation*}

Let $W_q$ and $\ff_q'$ be the same as the sets $W_r$ and $\ff_{r}'$ defined in the proof of the previous theorem, but now with $(\alpha_i):=(\alpha_i(q))$ and $n \ge \max \set{m_1, \ldots,m_k}$ large enough such that the inequality
\begin{equation}\label{24}
\left(1+\frac{1}{q^n}\right)^k < 1 + \frac{1}{q^{m_1 + \cdots + m_k}}
\end{equation}
holds. If $w_{i_1j_1}, \ldots, w_{i_kj_k}$ are $k$ words belonging to $W_q$ such that
$$i_1j_1 \ldots i_kj_k \not=\ell_1 m_1 \ldots \ell_km_k,$$
we associate with them a similarity $S_{i_1j_1 \ldots i_kj_k} : J_q \to J_q$ defined by the formula
\begin{align*}
S_{i_1j_1 \ldots i_kj_k}(x) & = \frac{\alpha_1}{q} + \cdots + \frac{\alpha_{j_1-1}}{q^{j_1-1}} + \frac{i_1}{q^{j_1}} \\
& + \frac{\alpha_1}{q^{j_1+1}} + \cdots + \frac{\alpha_{j_2-1}}{q^{j_1 + j_2 -1}} + \frac{i_2}{q^{j_1 +j_2}} \\
&  \, \, \vdots \\
& + \frac{\alpha_1}{q^{j_1+ \cdots + j_{k-1} + 1}} + \cdots + \frac{\alpha_{j_k-1}}{q^{j_1 + \cdots + j_k -1}} + \frac{i_k}{q^{j_1 + \cdots + j_k}}
\\
& + \frac{x}{q^{j_1 + \cdots + j_k}}, \quad x \in J_q.
\end{align*}
Let $\mathcal{G}_q'$ denote the set of those sequences belonging to $\ff_q'$ which do not contain the word $w$, and let
\begin{equation*}
\mathcal{G}_q:= \set{\sum_{i=1}^{\infty} \frac{c_i}{q^i} : (c_i) \in \mathcal{G}_q'}.
\end{equation*}

Since $(\alpha_i)=(\alpha_i(q))$, a sequence belonging to $\ff_q'$ is not necessarily the greedy expansion in base $q$ of a number $x \in J_q$, but this does not affect our proof. It is important, however, that any greedy expansion $(b_i) \not= \alpha_1^{\infty}$ in base $q$ can be written as $\alpha_1^{\ell} c_1c_2 \ldots$ for some $\ell \ge 0$ and some sequence $(c_i)$ belonging to $\ff_q'$. If $Y_w$ denotes the set of numbers  $x \in J_q$ for which the word $w$ does not occur in $(b_i(x,q))$ then the latter fact implies that the set $Y_w \setminus \set{\alpha_1/(q-1)}$ can be covered by countably many sets similar to $\mathcal{G}_q$.

It follows from the definition of $\mathcal{G}_q$ that
\begin{equation*}
\mathcal{G}_q \subset \bigcup S_{i_1j_1 \ldots i_kj_k}(\mathcal{G}_q)
\end{equation*}
where the union runs over all $i_1j_1 \ldots i_kj_k$ for which the similarity $S_{i_1j_1 \ldots i_kj_k}$ is defined above.
Hence
\begin{equation*}
\overline{\mathcal{G}_q} \subset \bigcup S_{i_1j_1 \ldots i_kj_k}(\overline{\mathcal{G}_q})
\end{equation*}
and thus $\mathcal{G}_q \subset \hh_q$ where $\hh_q$ is the (nonempty compact) invariant set of this system of similarities.
Let $\tilde{\alpha_i}:= \alpha_i$ for $1 \le i < n$ and $\tilde{\alpha_n}:=\alpha_{n}+1$. From Proposition 9.6 in \cite{F} we know that ${\rm dim}_{H} \hh_q \le s$ where $s$ is the real solution of the equation
\begin{equation}\label{25}
\sum_{j_1=1}^n \sum_{j_2=1}^n \cdots \sum_{j_k=1}^n \left(\frac{\Pi_{i=1}^k \tilde{\alpha_{j_i}}}{q^{(j_1 + \cdots + j_k)s}}\right)
-\frac{1}{q^{(m_1+ \cdots + m_k)s}}=1.
\end{equation}
Denoting the left side of \eqref{25} by $C(s)$, we have

\begin{equation*}
C(1) + \frac{1}{q^{m_1 + \cdots + m_k}}  = \left(\sum_{i=1}^n \frac{\tilde{\alpha_{i}}}{q^{i}} \right)^k < \left(1 + \frac{1}{q^n} \right)^k.
\end{equation*}
By \eqref{24} we have $C(1) < 1$, and thus ${\rm dim}_{H} Y_v \le {\rm dim}_{H} Y_w \le {\rm dim}_{H} \hh_q < 1$.

(ii) An argument analogous to the one used in the proof of Theorem~\ref{t27}(ii) shows that the set $Y_w$ (and thus $Y_v$) is of first category. Since $Y$ is a countable union of sets of the form $Y_v$ it follows that $Y$ is a null set of first category.
Let $r \in (1,q)$ and let $G_{r,q}$ be the set defined in Theorem~\ref{t27}. Due to Theorem~\ref{t27}(i) it is now sufficient to show that $G_{r,q} \subset Y$. To this end, choose a number $p \in (r,q)$, and let $n \in \NN$ be large enough such that the inequalities
\begin{equation*}
\alpha_1(r) \ldots \alpha_n(r) < \alpha_1(p) \ldots \alpha_n(p) < \alpha_1(q) \ldots \alpha_n(q)
\end{equation*}
hold. Note that such an integer $n$ exists by Proposition~\ref{p21}. From Propositions~\ref{p21} and ~\ref{p24} we conclude that the sequence $0\alpha_1 (p) \ldots \alpha_n(p)0^{\infty}$ equals $(b_i(y,q))$ for some $y \in [0,1)$ while the word $0\alpha_1(p) \ldots \alpha_n(p)$ does not occur in the greedy expansion in base $r$ of any number $x \in J_r$.
\end{proof}

\section{Proof of Theorem~\ref{t12}}\label{s3}

The following characterization of unique expansions readily follows from Proposition \ref{p24}.

\begin{proposition}\label{p31}
Fix $q>1$. A sequence $(c_i)$ of integers $c_i \in A_q$ is the unique expansion of some $x\in J_q$ if and
only if
\begin{equation*}
c_{n+1}c_{n+2}\ldots <\alpha_1(q)\alpha_2(q)\ldots \quad\text{whenever }c_n< \alpha_1(q)
\end{equation*}
and
\begin{equation*}
\overline{c_{n+1}c_{n+2}\ldots }< \alpha_1(q)\alpha_2(q)\ldots \quad\text{whenever }c_n>0.
\end{equation*}
\end{proposition}

In what follows we use the notations $(a_i(x,q))$,
$(b_i(x,q))$, $(\alpha_i(q))$ and $(\beta_i(q))$ as introduced in Section~\ref{s1}. If $x$
and $q$ are clear from the context, then we omit these arguments and we simply write $a_i$, $b_i$, $\alpha_i$ and
$\beta_i$. If two couples $(x,q)$ and $(x',q')$ are considered simultaneously, then we also write $a'_i$, $b'_i$,
$\alpha'_i$ and $\beta'_i$ instead of $a_i(x',q')$, $b_i(x',q')$, $\alpha_i(q')$ and $\beta_i(q')$.

\begin{lemma}\label{l32}
Given $(x,q)\in\bfj $, the following two conditions are equivalent:
\begin{align*}
&\overline{a_{n+1}a_{n+2}\ldots }\le \alpha_1\alpha_2\ldots  \quad\text{whenever }a_n>0;\\
&\overline{a_{n+1}a_{n+2}\ldots }\le \beta_1\beta_2\ldots  \quad\text{whenever }a_n>0.
\end{align*}
\end{lemma}

\begin{proof}
Since $(\alpha_i) \le (\beta_i)$, it suffices to show that if there exists a positive integer $n$ such that
\begin{equation*}
a_n>0\quad\text{and}\quad \overline{a_{n+1}a_{n+2}\ldots }>\alpha_1\alpha_2\ldots ,
\end{equation*}
then there exists also a positive integer $n'$ such that
\begin{equation*}
a_{n'}>0\quad\text{and}\quad \overline{a_{n'+1}a_{n'+2}\ldots }>\beta_1\beta_2\ldots .
\end{equation*}
If the greedy expansion $(\beta_i)$ is infinite, then $(\beta_i)=(\alpha_i)$ and we may choose $n'=n$. If
$(\beta_i)$ has a last nonzero digit $\beta_{\ell}$, then
$(\alpha_i)=(\alpha_1\ldots \alpha_{\ell})^{\infty}$ with $\alpha_1\ldots \alpha_{\ell-1} \alpha_{\ell}=\beta_1\ldots \beta_{\ell -1} \beta_{\ell}^-$ $(\beta_{\ell}^-:=\beta_{\ell} -1)$, and
thus $\alpha_{\ell}<\alpha_1$. Since we have
\begin{equation*}
\overline{a_{n+1}a_{n+2}\ldots }>(\alpha_1\ldots \alpha_{\ell})^{\infty}
\end{equation*}
by assumption, there exists a nonnegative integer $j$ satisfying
\begin{equation*}
\overline{a_{n+1}\ldots a_{n+j\ell}}=(\alpha_1\ldots \alpha_{\ell})^j
\quad\text{and}\quad
\overline{a_{n+j\ell+1}\ldots a_{n+(j+1)\ell}}>\alpha_1\ldots \alpha_{\ell}.
\end{equation*}
Putting $n':=n+j\ell $ it follows that
\begin{equation*}
a_{n'}>0\quad\text{and}\quad \overline{a_{n'+1}\ldots a_{n'+\ell }}\ge\beta_1\ldots \beta_{\ell}.
\end{equation*}
It follows from our assumption $\overline{a_{n+1}a_{n+2}\ldots }>\alpha_1\alpha_2\ldots$ that
$(\alpha_i)<\alpha_1^{\infty}$ and $(a_i) \not= \alpha_1^{\infty}$. It follows from Proposition~\ref{p22} that
$(a_i)$ has no tail equal to $\alpha_1^{\infty}$, so that $\overline{a_{n'+\ell+1}a_{n'+\ell+2}\ldots }>0^{\infty}$. We conclude that
\begin{equation*}
\overline{a_{n'+1}a_{n'+2}\ldots }>\beta_1\beta_2\ldots .\qedhere
\end{equation*}
\end{proof}

\remove{
\begin{remark}
If we allow the digit $q$ for integer bases $q$, then Proposition \ref{p31} and Lemma~\ref{l32} remain valid by changing the
definition of the conjugate to $\overline{b_i(x,q)}:=\lfloor q\rfloor -b_i(x,q)$ and $\overline{a_i(x,q)}:=\lfloor q\rfloor
-a_i(x,q)$.
\end{remark}
}
\begin{definition}
We say that $(x,q)\in\bfj $ belongs to the set $\bfv $ if one of the two equivalent conditions of the preceding lemma
is satisfied. Moreover, we define
\begin{equation*}
\vv_q:= \set{ x \in J_q: (x,q) \in \bfv}, \quad q>1.
\end{equation*}
\end{definition}

It follows from Proposition~\ref{p31} that $\bfu \subset\bfv \subset\bfj $.

\begin{proof}[Proof of Theorem~\ref{t12}] We need to prove that $\bfuu \cap \bfj = \bfv$.

First we show that $\bfv \subset\bfuu $. In order to do so, we introduce for each fixed $q>1$ the sets $\uu_q'$ and  $\vv_q'$, defined by
\begin{equation*}
\uu_q':= \set{(a_i(x,q)): x \in \uu_q} \quad \text{and} \quad \vv_q':= \set{(a_i(x,q)) : x \in \vv_q}.
\end{equation*}
Observe that $\uu_q'$ is simply the set of unique expansions in base $q$. It follows easily from Propositions \ref{p21}, \ref{p22} and \ref{p31} that $\uu_q' \subset \vv_q'$ for each $q >1$, and that $\vv_q' \subset \uu_r'$ for each $r > q$ such that
$\lceil q \rceil = \lceil r \rceil$. Since we also have $\uuuq = \vv_q = [0,1]$ if
$q>1$ is an integer, the result follows.

\medskip
Next we show that $\bfuu \cap\bfj \subset\bfv $. Since $\bfu \subset \bfv$ it is sufficient to prove that
if $(x,q)\in\bfj \setminus\bfv $,  then $(x',q')\notin\bfv $ for  all $(x',q')\in\bfj $ close enough to $(x,q)$.
Applying Lemma \ref{l32} there exist two positive integers $n$ and $m$ such that
\begin{equation}\label{31}
a_{n}>0\quad\text{and}\quad \overline{a_{n+1}\ldots a_{n+m}}>\beta_1\ldots\beta_m .
\end{equation}
This implies in particular that $q$ is not an integer, because otherwise $(\alpha_i) = (\beta_i) = \beta_1^{\infty}$.
Hence, if $q'$ is sufficiently close to $q$, then
\begin{equation}\label{32}
\beta'_1\ldots \beta'_m\le \beta_1\ldots \beta_m
\end{equation}
by Lemma \ref{l25}.
It follows from the definition of quasi-greedy expansions that
\begin{equation*}
\frac{a_1}{q}+\cdots+\frac{a_{j-1}}{q^{j-1}}+ \frac{a_j^+}{q^j}+\frac{1}{q^{j+m}}>x\quad\text{whenever $a_j<\alpha_1$,}
\end{equation*}
where $a_j^+:=a_j+1$. If $(x',q') \in \bfj$ is sufficiently close to $(x,q)$, then $\alpha_1=\alpha_1'$, the inequality \eqref{32} is satisfied, $a_1' \ldots a_{n+m}' \ge a_1 \ldots a_{n+m}$ by Lemma \ref{l23}(i),
and
\begin{equation}\label{33}
\frac{a_1}{q'}+\cdots+\frac{a_{j-1}}{(q')^{j-1}}+ \frac{a_j^+}{(q')^j}+\frac{1}{(q')^{j+m}}>x'\quad\text{whenever $j\le n+m$ and
$a_j<\alpha_1$.}
\end{equation}
Now we distinguish between two cases.

If $a_1'\ldots a'_{n+m} =  a_1\ldots a_{n+m}$, then we have
\begin{equation*}
a'_n>0\quad\text{and}\quad \overline{a'_{n+1}\ldots a'_{n+m}}>\beta_1\ldots\beta_m\ge  \beta'_1\ldots \beta'_m
\end{equation*}
by \eqref{31} and \eqref{32}. This proves that $(x',q')\notin\bfv $.

If $a_1'\ldots a'_{n+m}> a_1\ldots a_{n+m}$, then let us consider the smallest $j$ for which $a'_j>a_j$. It follows from \eqref{32} and \eqref{33} that
\begin{equation*}
a'_j=a_j^+>0\quad\text{and}\quad \overline{a'_{j+1}\ldots a'_{j+m}}=\beta_1^m>\beta_1\ldots\beta_m\ge  \beta'_1\ldots \beta'_m.
\end{equation*}
Hence $(x',q')\notin\bfv $ again.
\end{proof}

\begin{remark}
It is the purpose of this remark to describe the set $\bfuu \setminus \bfj$. For each $m \in \NN$, we define the number $q_m \in (m, m+1)$ by the equation
\begin{equation*}
1 = \frac{m}{q_m} + \frac{1}{q_m^2}.
\end{equation*}
Fix $q \in (m, q_m]$.
Since $\alpha_1(q)=m$ and $\alpha_2(q)=0$, Proposition~\ref{p31} implies that a sequence $(c_i) \in \set{0, \ldots, m}^{\NN}$ belongs to $\uu_q'$ if and only if for each $n \in \NN$, we have
\begin{equation*}
c_n < m \Longrightarrow c_{n+1} < m
\end{equation*}
and
\begin{equation*}
c_n >0 \Longrightarrow c_{n+1} >0.
\end{equation*}
Denoting the set of all such sequences by $D_m'$ and putting for $m >1$ (note that $D_1'=\set{0^{\infty},1^{\infty}}$),
\begin{equation*}
D_m:= \set{\sum_{i=1}^{\infty} \frac{c_i}{m^i}: (c_i) \in D_m'},
\end{equation*}
one may verify that
\begin{equation*}
\bfuu \setminus \bfj=\set{(0,1)}\cup \bigcup_{m=2}^{\infty} (D_m \setminus [0,1]) \times \set{m}.
\end{equation*}
\end{remark}

\section{Proof of Theorem \ref{t11}}\label{s4}

We need some results on the Hausdorff dimension of the sets $\uu_q$ and $\vv_q$ for $q>1$.
It follows from Theorem~\ref{t12} that $\uu_q \subset \uuuq \subset \vv_q$. Moreover, if an element $x \in \vv_q \setminus \uu_q$ has an infinite greedy expansion in base $q$, then $(b_i(x,q))$ must end with $\overline{\alpha_1(q)\alpha_2(q) \ldots}$ as follows from Propositions~\ref{p24} and \ref{p31}; hence $\vv_q \setminus \uu_q$ is countable and the sets $\uu_q$, $\uuuq$ and $\vv_q$ have the same Hausdorff dimension for each $q>1$. Proposition~\ref{p41} below is contained in the works of Dar\'oczy and K\'atai \cite{DK2}, Kall\'os \cite{K1}, \cite{K2}, Glendinning and Sidorov  \cite{GS}, and Sidorov \cite{Si2}; for the reader's convenience we provide here an elementary proof.

\begin{proposition}\label{p41} We have
\begin{itemize}
\item[\rm (i)] $\lim_{q\uparrow 2}{\rm dim}_{H} \uu_q=1$;
\item[\rm (ii)] ${\rm dim}_{H} \uu_q < 1$ for all non-integer $q >1$.
\end{itemize}
\end{proposition}

\begin{proof}
(i) Assume that $q \in (1,2)$ is larger than the tribonacci number, i.e.,
\begin{equation*}
\frac{1}{q}+\frac{1}{q^2}+\frac{1}{q^3}<1,
\end{equation*}
and let $N=N(q) \ge 2$ be the largest integer satisfying
\begin{equation*}
\frac{1}{q}+\cdots+\frac{1}{q^{2N-1}}<1.
\end{equation*}
Hence $\alpha_1(q)=\cdots=\alpha_{2N-1}(q)=1$.
Let us denote by $\mathcal{I}_q$ the set of numbers $x\in J_q$ which have an expansion $(c_i)$ satisfying $0<c_{kN+1}+\cdots
+c_{kN+N}<N$ for every nonnegative integer $k$. Since in such expansions $(c_i)$, a zero (one) is followed by at most
$2N-2$ consecutive one (zero) digits, it follows from Proposition~\ref{p31} that
$\mathcal{I}_q \subset \uu_q$. Moreover, the set $\mathcal{I}_q$ is closed and thus compact. In order to prove this, observe first that the set $\mathcal{I}_q$ is closed from above~\footnote{We call a set $X \subset \RR$ \emph{closed from above} if $x$ belongs to $X$ whenever there exists a sequence $(x_n)$ of elements of $X$ that converges to $x$ from above.} by virtue of Lemma~\ref{l25}(ii). Hence $\mathcal{I}_q$ is closed because $\mathcal{I}_q$ is symmetric relative to $J_q$.

It suffices now to prove that
\begin{equation}\label{41}
{\rm dim}_{H} \mathcal{I}_q= \frac{\log (2^N-2)}{N\log q};
\end{equation}
indeed, $q\uparrow 2$ implies $N\to\infty$, hence ${\rm dim}_{H} \mathcal{I}_q\to 1$ and consequently ${\rm dim}_{H} \uu_q\to 1$.

Observe that
\begin{equation}\label{42}
\mathcal{I}_q=\bigcup S_{c_1\ldots c_N}(\mathcal{I}_q)
\end{equation}
where the union runs over the words $c_1 \ldots c_N$ of length $N$ consisting of zeros and ones satisfying $0<c_1+\cdots +c_N<N$, and $S_{c_1 \ldots c_N} : J_q \to J_q$ is given by
\begin{equation*}
S_{c_1\ldots c_N}(x):=\Bigl(\frac{c_1}{q}+\cdots+\frac{c_N}{q^N}\Bigr)+\frac{x}{q^N}, \quad x \in J_q.
\end{equation*}
In other words, $\mathcal{I}_q$ is the (nonempty compact) invariant set of the iterated function system formed by these $2^N-2$ similarities. The
sets $S_{c_1\ldots c_N}(\mathcal{I}_q)$ on the right side of \eqref{42} are disjoint because $S_{c_1 \ldots c_N}(\mathcal{I}_q) \subset \mathcal{I}_q \subset \uu_q$, and
since all similarity ratios are equal to $q^{-N}$, it follows from Propositions 9.6 and 9.7 in \cite{F}
that the Hausdorff dimension $s$ of $\mathcal{I}_q$ is the real solution of the equation
\begin{equation*}
(2^N-2)q^{-Ns}=1,
\end{equation*}
which is equivalent to \eqref{41}.

(ii) Let $q>1$ be a non-integer and let $n \in \NN$ be such that $\alpha_n(q) < \alpha_1(q)$. It follows from Proposition~\ref{p31} that the word $1(0)^n$ does not occur in $(b_i(x,q))$ if $x$ belongs to $\uu_q$. Applying Theorem~\ref{t28}(i) with $y=q^{-1}$, $\ell=0$ and $m=n+1$, we conclude that ${\rm dim}_{H} \uu_q < 1$.
\end{proof}

\begin{proof}[Proof of Theorem \ref{t11}]
(ii) Let $q>1$ be a non-integer. Since $\vv_q \setminus \uu_q$ is countable, Proposition~\ref{p41}(ii) yields that ${\rm dim}_{H} \vv_q <1$. This implies in particular that the set $\vv_q$ is a one-dimensional null set. Applying Theorem~\ref{t12} (and the remark following its proof) and Fubini's theorem we conclude that $\bfuu$ is a two-dimensional null set.

(i)  Since $\uu_q$ is not closed for all $q>1$, $\bfu$ cannot be closed. Since $\bfuu$ is a two-dimensional null set, it has no
interior points. It remains to show that $\bfu$ (and thus $\bfuu$) has no isolated points. If $q > 1$ is an integer, then, as is well known, $\uu_q$ is dense in $J_q=[0,1]$. If $q>1$ is a non-integer, then each $(x,q) \in \bfu$ is not isolated because $\uu_q' \subset \uu_r'$ whenever $q < r$ and $\lceil q \rceil = \lceil r \rceil$.

(iii) From Corollary 7.10  in \cite{F} we may conclude that for almost all $q>1$
\begin{equation*}
{\rm dim}_H \uu_q\le \max\set{0,{\rm dim}_H \bfu  -1}
\end{equation*}
which would contradict Proposition~\ref{p41}(i) if we had ${\rm dim}_H \bfu< 2$.
\end{proof}

\noindent
{\it Acknowledgements.} The first author has been supported by NWO, Project nr. ISK04G.
Part of this work was done during a visit of the second author at the Department of Mathematics of the Delft
University of Technology. He is grateful for this invitation and for the excellent working conditions.


\begin{thebibliography}{[K$^3$1]}

\bibitem{BK} C. Baiocchi, V. Komornik,
\emph{Greedy and quasi-greedy expansions in non-integer bases},
arXiv:0710.3001 [math.], October 16, 2007.

\bibitem{DDV} K. Dajani, M. de Vries,
\emph{Invariant densities for random $\beta$-expansions},
J. Eur. Math. Soc. {\bf 9} (2007), no. 1, 157--176.

\bibitem{DK1} Z. Dar\'oczy, I. K\'atai,
\emph{Univoque sequences},
Publ.\ Math.\ Debrecen {\bf 42} (1993), no. 3--4, 397--407.

\bibitem{DK2} Z. Dar\'oczy, I. K\'atai,
\emph{On the structure of univoque numbers},
Publ.\ Math.\ Debrecen {\bf 46} (1995), no. 3--4, 385--408.

\bibitem{DVK} M. de Vries, V. Komornik,
\emph{Unique expansions of real numbers},
Adv. Math. {\bf 221} (2009), no. 2, 390--427.

\bibitem{EHJ} P. Erd\H{o}s, Horv\'ath, M., I. Jo\'o,
\emph{On the uniqueness of the expansions $1=\sum q^{-n_i}$},
Acta Math. Hungar. {\bf 58} (1991), no. 3--4, 333--342.

\bibitem{EJK1} P. Erd\H{o}s, I. Jo\'o, V. Komornik,
\emph{Characterization of the unique expansions} $1=\sum\sp \infty\sb {i=1}q\sp {-n\sb i}$ \emph{and related problems},
Bull.\ Soc.\ Math.\ France {\bf 118} (1990), no. 3, 377--390.

\bibitem{EJK2} P. Erd\H{o}s, I. Jo\'o, V. Komornik,
\emph{On the number of $q$-expansions},
Ann.\ Univ.\ Sci.\ Budapest.\ E\"otv\"os Sect.\ Math.\ {\bf 37} (1994), 109--118.

\bibitem{F} K. Falconer,
\emph{Fractal Geometry. Mathematical Foundations and Applications},
John Wiley \& Sons, Chicester, second edition, 2003.

\bibitem{GS} P. Glendinning, N. Sidorov,
\emph{Unique representations of real numbers in non-integer bases},
Math.\ Res.\ Lett. {\bf 8} (2001), no. 4, 535--543.

\bibitem{K1} G. Kall\'os,
\emph{The structure of the univoque set in the small case},
Publ.\ Math.\ Debrecen {\bf 54} (1999), no. 1--2, 153--164.

\bibitem{K2} G. Kall\'os,
\emph{The structure of the univoque set in the big case},
Publ.\ Math.\ Debrecen {\bf 59} (2001), no. 3--4, 471--489.

\bibitem{KK} I. K\'atai, G. Kall\'os,
\emph{On the set for which $1$ is univoque},
Publ. Math. Debrecen {\bf 58} (2001), no. 4, 743--750.

\bibitem{KL1} V. Komornik, P. Loreti,
\emph{Unique developments in non-integer bases},
Amer. Math. Monthly {\bf 105} (1998), no. 7, 636--639.

\bibitem{KL2} V. Komornik, P. Loreti,
\emph{On the topological structure of univoque sets},
J. Number Theory, {\bf 122} (2007), no. 1, 157--183.

\bibitem{P} W. Parry,
\emph{On the $\beta $-expansions of real numbers},
Acta Math. Acad. Sci. Hungar. {\bf 11} (1960), 401--416.

\bibitem{R} A. R\'{e}nyi,
\emph{Representations for real numbers and their ergodic properties},
Acta Math. Acad. Sci.\ Hungar. {\bf 8} (1957), 477-493.

\bibitem{Si1} N. Sidorov,
\emph{Almost every number has a continuum of $\beta$-expansions},
Amer. Math. Monthly {\bf 110} (2003), no. 9, 838--842.

\bibitem{Si2} N. Sidorov, \emph{Combinatorics of linear iterated function systems with overlaps}, Nonlinearity {\bf 20} (2007), 1299--1312.
\end{thebibliography}
\end{document}